\documentclass{preprint}

\usepackage{enumerate} 				% roman enumeration and hyperref
\usepackage{hyperref} 				% enables hyperlinks
\usepackage{tikz-cd}				% for diagrams
\usepackage{frege}				    % for typesetting symbols in the ideography-style

\title{Sense, reference, and computation}

\author{Bruno Bentzen}

\department{Department of Philosophy}

\university{Carnegie Mellon University}

\country{Pittsburgh, PA \\ USA}

\email{b.bentzen@hotmail.com}

\journal{Perspectiva Filos\'ofica}

\volume{47}

\issue{2}

\when{2020}
\pages{179--203}
\setdoi{}

\defcitealias{hottbook}{UFP, 2013}

\begin{document}

\nocite{frege1879begriffsschrift,frege1884grundlagen,frege1892sinn,frege1962grundgesetze}

\maketitle

\begin{abstract}
In this paper, I revisit Frege's theory of sense and reference in the constructive setting of the meaning explanations of type theory, extending and sharpening a program--value analysis of sense and reference proposed by Martin-L\"of building on previous work of Dummett. I propose a computational identity criterion for senses and argue that it validates what I see as the most plausible interpretation of Frege's equipollence principle for both sentences and singular terms. Before doing so, I examine Frege's implementation of his theory of sense and reference in the logical framework of \gga, his doctrine of truth values, and views on sameness of sense as equipollence of assertions. 
\end{abstract}

%=========================================================================================================================================
\section{Introduction} \label{intro}
%=========================================================================================================================================

Frege was one of the first to dispute that all mathematical truths are based on intuition, the cornerstone of the dominant Kantian philosophy of mathematics, and put forward the thesis that all arithmetical truths are reducible to logical truths.\footnote{
	References to Frege's works will be cited by name only and indicated by page or section number in the text. I shall ambiguously refer to the formal language and formal theory of \gga~by ``ideography'', unless explicitly stated otherwise.}  
The success of his logicism would eventually depend on the plausibility of the development of a new system of logic in which all the objects and concepts of arithmetic could be defined. Indeed, that was the original purpose of his \bs, the book which lays down the first systematic treatment of modern logic, the ideography, and applies the system to the mathematical theory of sequences. Frege hints at his logicist program in the preface to his book, announcing that his investigations will be continued with the logical elucidation of the concepts of number, magnitude, and so forth in what he describes as an immediately following publication. Yet, it would not be until fourteen years later that the promised book would finally appear in the form of the first volume of \gga. Frege himself explains that this long delay was due to the accommodation of essentially two major changes to the ideography, a task that ultimately forced him to discard an almost completed manuscript.\footnote{Cf. \gga~I, ix.} The first change was the introduction of value-ranges of functions, new objects that are intended to generalize the informal notion of extension of concept that he invokes in {\S68} of \gla~to overcome the infamous Julius Caesar problem with an explicit definition of the concept of number. This point has been studied extensively and it will not concern us here.\footnote{
	For a recent survey on the topic, see \cite{bentzen2019frege}.} 
Instead, I would like to focus on Frege's theory of sense and reference, the doctrine that singular terms and complete sentences have a mode of presentation and denote an object, whose implementation consists in the second substantial change made to the ideography. 

Unfortunately, Frege's formal account of his distinction between sense and reference leaves much to be desired, as he failed to develop his ideography sufficiently in order to capture his answer to the paradox of identity, namely, the problem of explaining the apparent difference in cognitive value between true identity statements of the form $a = a$ and $a = b$. This will be the main subject of \sectionref{frege-sense}, where I shall examine the role of the theory of sense and reference in the technical development of Frege's logicist project and his views on sameness of sense as equipollence of assertions. In \sectionref{intuitionistic-semantics}, I shall shift gears and turn to an opposing philosophical position that I take to possess a more suitable semantic setting for the establishment of a proper theory of sense and reference, namely, the constructive conception of mathematics. The idea of reinterpreting Frege's traditional semantic distinction in a constructive setting is certainly not new, as it has been first explored mainly by \cite{dummett1978sense} and \cite{martinlof2001sense}. The latter author has extended considerably the ideas of the first one by reshaping them in the context of his type theory, a formal system in which every term is assigned to a type and every operation strictly restricted to terms of a certain type.\footnote{
	In this paper, I will be concerned with a fragment of \cite{martinlof1982constructive} where the only types are dependent function types, equality types, and the natural numbers. Judgments of typehood are used instead of type universes, and typehood and type membership judgments are defined via primitive type and term judgments, as described in \sectionref{intuitionistic-semantics}.} 
But at the same time several aspects of Martin-L\"of's account strike me as incomplete or problematic, especially because he does not seem to make full use of the powerful intended semantics of type theory, an informal realizability interpretation known as the meaning explanations~\citep{martinlof1982constructive}. My goal in this paper is to revise some of Martin-L\"of's main views, a task that will be carried out in \sectionref{sense-and-reference}, where I carry the meaning explanations to their ultimate logical conclusions. 

%=========================================================================================================================================
\section{Frege on sameness of sense} \label{frege-sense}
%=========================================================================================================================================

In Frege's works we find the two forms of judgment that are current in modern logic, the judgment form that asserts that a proposition is true, which is  written in turnstile notation as
$$\Fa[1] A $$
and the judgment form that states that a sentence expresses a proposition, which is implicitly introduced for the first time in \S2 of \bs~via the content stroke 

$$\F[1] A.$$

\noindent Therein, a proposition is described as a ``judgeable content'' and identified with the content of a turnstile judgment. This intuitive notion of content represents the main semantic element of the ideography at this earlier point and, in fact, logical operations such as implication, negation, and universal quantification are all exclusively applicable to contents, or, more precisely, to judgeable contents or propositions. 

The only exception to this structure is the relation of identity of content which is introduced in \S8 of \bs~as a metalinguistic relation between the expressions themselves and not their respective contents. That is, according to this interpretation of equality, an equality statement $a = b$ means that the expressions `$a$' and `$b$' have the same content. One semantic consequence of this is that two propositions cannot be identified if they imply each other, for they may have different contents. Indeed, the propositional extensionality principle 

$$ (A \varsupset B) \varsupset (B \varsupset A) \varsupset (A = B)$$

\noindent is not only not a theorem in \bs~but also inconsistent with its formal system, as noted in \cite{duarte2009logica}. However, once the ambiguous notion of content is split into sense and reference, propositional extensionality becomes a valid principle for references, although not for senses, considering that from the fact that two theorems are logically equivalent it does not follow that they must have the same cognitive value.\footnote{
	Curiously, \cite{frege1980correspondence} once considered an identity criterion for senses along those lines in a letter to Husserl dated 9 December 1906, but this suggestion is clearly unacceptable, since it would mean that every two true propositions $a = a$ would $a = b$ express the same sense.} 

%-----------------------------------------------------------------------------------------------------------------------------------------
\subsection{The doctrine of truth values}
%-----------------------------------------------------------------------------------------------------------------------------------------  

Frege's theory of sense and reference is systematically introduced for both singular terms and complete sentences in his seminal essay \textit{Sinn und Bedeutung}. The theory seems to provide a compelling answer to the paradox of identity with respect to singular terms, for we are able to say that the cognitive value of $a = a$ and $a = b$ is generally different because $a$ and $b$ may differ in their sense, despite being coreferential. However, the theory becomes highly controversial when we consider the sense and reference of sentences, which are taken to express a thought (in modern terminology, a proposition) and refer to a truth value, which, since referents are always objects, are seen as the simplest representatives of logical objects. Frege is well aware of this counterintuitive aspect of his doctrine, as he immediately attempts to justify it by noting that what he calls an object can only be properly understood in connection to concepts and relations and then attempts to substantiate his claim that the reference of a sentence is its truth value by remarking that the latter remain unchanged when a part of the sentence is replaced by a coreferential expression, that is, that sentences with the same truth value are interchangeable \textit{salva veritate}. Later on, in the preface of \gga, Frege remarks that the introduction of truth-values may seem strange at first, but adds that the fact that everything becomes much simpler and sharper with those objects puts a great weight in the balance in favor of his own conception.\footnote{
	Cf. \gga~I, x.} 
In the ideography, concepts and relations are treated as particular cases of functions that output truth values and Frege explicitly states that an object is anything that is not a function.

Besides the poor justifications for Frege's adoption of his strange doctrine of truth values, it is disappointing that almost no considerations about sense or thoughts can be found in \gga, at least if we assume that the real purpose of the distinction between sense and reference is to offer a solution to the paradox of identity. Although Frege still seems to have a strong conviction that a criterion of identity for senses is of fundamental importance, as we can judge from, for instance, a December 1906 letter to Husserl, where Frege writes: 

\begin{quote}
	It seems to me that an objective criterion is necessary for recognizing a thought again as the same, for without it logical analysis is impossible. \cite[p. 70]{frege1980correspondence}
\end{quote}

\noindent Frege does not seem to bother to provide a logic of sense and reference in the book he intends to show once and for all that arithmetic is nothing but logic. For this reason, we can say that Frege failed to show that some mathematical equality statements have cognitive value in \gga, because his solution to the identity paradox is entirely dependent on what it means for two expressions to have different senses, and no identity criterion for senses is given in his ideography.

It is possible that Frege had come to recognize that the most important notion for the vindication of his logicism is that of reference, as pointed out by \cite{simons1992there}, and that for those purposes all that we need to know about sense is that a sentence needs to express a thought in order to refer at all. 
In fact, now the horizontal stroke symbolism 

$$\F[1] A$$

\noindent which used to indicate that $A$ is a proposition (judgeable content) before the division of content into sense and reference, no longer says that $A$ expresses a proposition (thought). Instead of being an implicit form of judgment, the horizontal is treated as a function that refers to a truth value depending on the reference of $A$, yielding the true if $A$ denotes the true and the false otherwise. 
Similarly, negation does not operate on contents anymore, for it is explicitly treated as a function 

$$\Fn[1] A$$
that has the opposite effect of yielding the false if $A$ does not refer to the true and the true otherwise. The turnstile judgment has kept its role, being the only judgment form of the revised ideography. It is now fully explained in terms of reference as well, now asserting that the expression $A$ refers to the true. There is no way to formally assert the fact that a sentence expresses a thought in \gga, which is to say to fulfill the role previously played by the content stroke. 

In general, it would seem that Frege contradicts himself by downplaying the value of sense and treating his doctrine of truth values as a mere technical device for the successful accomplishment of his logicism, as argued in \cite{ruffino1997wahrheitswerte} and \cite{duarte2009logica}. This can be seen more clearly in light of the discarded manuscript that was written before the introduction of the distinction between sense and reference to the ideography in \gga. 
In~\S69~of~\gla, immediately after the explicit definition of the concept of number is proposed by Frege, he attempts to confirm the fruitfulness of his definition by sketching a proof of a crucial theorem that says that the number that falls under the concept $F$ is equal to the number that falls under the concept $G$ iff $F$ and $G$ are in one-to-one correspondence (Hume's Principle). 
The proof consists in showing that both sides of the biconditional imply each other. But in the ideography, where a biconditional is always represented by an equality, the formalization of the proof cannot be completed without propositional extensionality. In \gga, what makes the derivation of propositional extensionality as a theorem possible is the introduction of Basic Law IV, which roughly states that the truth values denoted by the sentences $A$ and $B$ either coincide or not:
$$\neg (\F[1] A = \Fn[1] B) \varsupset (\F[1] A = \F[1] B).$$ 
But this is only allowed as an axiom in the ideography because of the distinction between sense and reference. If those observations are correct, Frege conceived his theory of sense and reference primarily as a logical apparatus necessary to overcome technical obstacles rather than as a solution to the paradox of identity as it is commonly thought.\footnote{
	See e.g. \citet[\S3]{ruffino1997wahrheitswerte} and \citet[p.167]{duarte2009logica}.} 
That would explain why Frege does not do justice to the notion of sense in \gga. 

%-----------------------------------------------------------------------------------------------------------------------------------------
\subsection{The equipollence principle}
%-----------------------------------------------------------------------------------------------------------------------------------------

What comes closer to the proposal of a criterion of identity for senses in Frege's writings is the equipollence principle that two sentences $A$ and $B$ express the same thought provided that anyone who accepts $A$ as true must also immediately accept $B$ as true and vice versa. \cite{sundholm1994frege} has called attention to the fact that this idea has been informally suggested on many occasions by Frege, most notably in the opening passage of his \textit{Kurze \"Ubersicht meiner logischen Lehren}: 

\begin{quote}
	Now two propositions $A$ and $B$ can stand in such a relation that anyone who recognizes the content of $A$ as true must thereby also recognize the content of $B$ as true and, conversely [...] So one has to separate off from the content of a proposition the part that alone can be accepted as true or rejected as false. I call this part the thought expressed by the proposition.
	\cite[pp.197--98]{frege1906brief}
\end{quote}

\noindent I take the following equipollence to be the most natural choice of a corresponding principle for singular terms: two singular terms $a$ and $b$ express the same sense if for every predicate $P$, anyone who accepts $P(a)$ as true must also immediately accept $P(b)$ as true and vice versa. Related to \textit{salva veritate}, this idea is already hinted at in \textit{Sinn und Bedeutung}

\begin{quote}
	If we now replace one word of the sentence by another having the same reference, but a different sense, this can have no bearing upon the reference of the sentence. Yet we can see that in such a case the thought changes; since, e.g., the thought in the sentence `The morning star is a body illuminated by the Sun' differs from that in the sentence `The evening star is a body illuminated by the Sun.' Anybody who did not know that the evening star is the morning star might hold the one thought to be true, the other false. \citep[p.62]{frege1892sinn}
\end{quote}

Yet, much remains to be done before both equipollence principles for singular terms and sentences can be adopted as satisfactory identity criteria for senses, for we do not have a rigorous account of what it means to ``immediately'' accept a proposition as true. In \textit{Funktion und Begriff}, Frege has expressly stated that the sense of two expressions is equal up to renaming of bound variables, that is, that we could write $x^2 - 4x$ directly as $y^2 - 4y$ without altering its sense, however, it can be very difficult to determine the extent of strictness in his account of sameness of sense, considering that Frege also happens to claim that the two halves of Basic Law V express the same sense, but in a different way.\footnote{
	Both claims are found in~\citet[p.27]{frege1891funktion}. Frege's views on sameness of sense have been dealt with in more detail in  \citet[p.304--307]{sundholm1994frege}, \citet[\S4]{klement2016grundgesetze}, and \citet[\S2.1]{bentzen202Xdifferent}.} 

%=========================================================================================================================================
\section{Constructive semantics in type theory} \label{intuitionistic-semantics} 
%=========================================================================================================================================

In the remainder of this paper, I will argue that a constructive rendering of the theory of sense and reference is capable of not only providing a precise computational interpretation to the equipollence principles for singular terms and sentences, but also validating them. 
I will begin our discussion with a brief outline of some important aspects of constructive semantics. 

In Frege's semantics, as mentioned in the previous section, the assertion that a proposition is true is understood as declaring the fact that the thought expressed by a sentence denotes the true, an interpretation that rests on a realist assumption that numbers and other sorts of mathematical objects are non-physical and non-mental entities that do not exist in space or time. Therefore, the most fundamental form of judgment in logic and mathematics is the truth judgment, which Frege writes as:

$$ \Fa[1] A. $$

\noindent In contrast, the constructivist tradition states that the only legitimate way of understanding a proposition is as a specification of a construction with certain given properties and, moreover, that in order to assert that a proposition is true one has to exhibit a construction that realizes the specification expressed by the given proposition~\citep{martinlof1996meanings}. The central form of judgment is therefore the one that states that $a$ is a construction that realizes a proposition $A$, which, using the language of type theory, is often expressed as:

$$ a : A.$$

Put differently, to assert that a proposition is realized by a construction is to declare the truth of that proposition. One thing that should be emphasized is that the existence of a construction should not be understood as an ontological existential claim that may be verified by means independent of our knowledge, but rather as the epistemic act of conceiving such a mathematical construction. Accordingly, the constructive conception of truth can be regarded as a form of anti-realism with respect to mathematical objects. 

%-----------------------------------------------------------------------------------------------------------------------------------------
\subsection{Judgments}
%-----------------------------------------------------------------------------------------------------------------------------------------

In type theory, a construction is a humanly computable procedure or program that is generally accepted as a method for obtaining a mathematical object. Following the so-called propositions-as-types correspondence~\citep{howard1980formulae}, both propositions and sets are specifications of constructions that are structured through the unifying concept of a type. The elements of a type are commonly called terms, and, as expected, terms are interpreted as constructions. 

For a more concise semantic explanation, it is convenient to base our type theory on two equality forms of judgment that state that two expressions are equal types and also that two expressions are equal terms of a type. I shall write those equalities with a triple bar notation to emphasize the fact that they are judgmental relations, meaning that they only occur in assertions and therefore are not subject to operations such as those determined by logical connectives:

$$A \equiv B \type \qquad\qquad a \equiv b : A.$$

\noindent The reflexivity of those equality judgments can be used to define their more common counterparts $A \type$ and $a : A$, which say that an expression is a type and that an expression is a term of a type, as $A \equiv A \type$ and $a \equiv a : A $. In other words, to determine their meaning, it suffices to determine the meaning of the type and term equality as primitive judgments. 
The general strategy is that the meaning of a judgment is given via an untyped model of computation, an idea that derives directly from the meaning explanations~\citep{martinlof1982constructive}. Every meaningful statement is based on a notion of computation that is accepted as primitive and used to give terms a computational behavior, define types as term specifications, and to assign terms to types based on the values to which they compute. 

%-----------------------------------------------------------------------------------------------------------------------------------------
\subsection{Computation}
%-----------------------------------------------------------------------------------------------------------------------------------------

The whole process of meaning explanation starts with the specification of a native computation system that takes the form of a programming language, typically the untyped lambda-calculus extended with constants for the type formers, constructors, and eliminators of the type theory. In this paper, I will consider dependent function types $\Pi_{x:A}{B}$, equality types $a =_A b$, and the natural numbers $\nat$. First, we specify the syntax of the programming language: 

\begin{align*}
\mathsf{var}\,\,\, :=\: & x \:\vert\: y \:\vert\: z \:\vert\: ... \:\vert\: x' \:\vert\: y' \:\vert\: z' \:\vert\: ... \\
\mathsf{expr} :=\: & \mathsf{var} \\
	& \pitype{\mathsf{var} : \mathsf{expr}} \mathsf{expr} \:\vert\: \fun{\mathsf{var}}\mathsf{expr} \:\vert\: \mathsf{expr}(\mathsf{expr})  \\
	& \mathsf{expr} =_{\mathsf{expr}} \mathsf{expr} \:\vert\:  \refl(\mathsf{expr}) \:\vert\: \eqrec(\mathsf{expr},\mathsf{expr},\mathsf{expr}) \\
	& \nat \:\vert\: \zero \:\vert\: \succ(\mathsf{expr}) \:\vert\: \natrec(\mathsf{expr},\mathsf{expr},\mathsf{expr}). 
\end{align*}

\vspace{4mm}

\noindent Then, we endow this language with an operational semantics, a complete description of how terms are expected to compute given in the form of a transition relation over closed terms, that is, expressions with no occurrence of free variables, and reflexive on certain closed terms, which are regarded as execution values.

This transition relation is fully captured by the computation rules that I now symbolically describe, where $\shouldsf{a} \longmapsto \shouldsf{a'} $ indicates the fact that $\shouldsf{a}$ transitions to $\shouldsf{a'}$ and $\shouldsf{a} \val$ means that $a$ is a value, which is to say that $\shouldsf{a} \longmapsto \shouldsf{a}$ is the case. I will start considering the rules for the dependent function type, which generalizes the notion of universal quantification and indexed product of sets:

$$ \frac{}{\prod_{x : \shouldsf{A}} \shouldsf{B} \val} \quad \frac{}{\fun{x} \shouldsf{a} \val} $$

$$ \frac{\shouldsf{a} \longmapsto \shouldsf{a'}}{\shouldsf{a(b)} \longmapsto \shouldsf{a'(b)}} $$
$$\frac{}{(\fun{x}\shouldsf{a})(\shouldsf{b}) \longmapsto \shouldsf{a}[\shouldsf{b}/x]} \quad \frac{}{\fun{x} \shouldsf{a}(x) \longmapsto \shouldsf{a}.} $$

\vspace{4mm}

\noindent The dependent function type $\Pi_{(x : \shouldsf{A})} B$ is inhabited by functions $\fun{x}a$, where $a$ is an open term that may depend on $x$, and if $f$ is a function and $a$ a term, then $f(a)$ is the application of $f$ to $a$. The last two computation rules are of particular importance, because they induce some obvious resemblances between lambda terms $\fun x f(x)$ and Frege's value-ranges $\acute{\epsilon}f(\epsilon)$, since $\acute{\epsilon}f(\epsilon) \cap a = f(a)$ is a theorem in the ideography (see \gga~\S34), although $\acute{\epsilon}f(\epsilon) = f$ contradicts the intended semantics of the ideography since value-ranges are taken to be objects and cannot be functions. I explore this connection in detail in \cite{bentzen2020frege}, which contains a more comprehensive type-theoretic study of the ideography. 

The following rules provide a full computational account of the equality type $a =_{A} b$, the type-theoretic counterpart of the usual notion of an equality proposition. Unlike the judgmental equality $a \equiv b : A$, $a =_A b$ does not have any assertive force. The constructor $\refl(a)$ represents the reflexivity of equality and the eliminator $\eqrec(a,b)$ generalizes the principle of \textit{salva veritate} that Frege borrows from Leibniz:

$$ \frac{}{\shouldsf{a =_{A} b} \val} \quad \frac{}{\refl(\shouldsf{a}) \val} $$
$$\frac{\shouldsf{a \longmapsto a'}}{\eqrec(\shouldsf{a,b}) \longmapsto \eqrec(\shouldsf{a',b})}$$
$$ \frac{}{\eqrec(\refl(\shouldsf{a}),\shouldsf{b}) \longmapsto \shouldsf{b}.}$$

\vspace{4mm}

\noindent Finally, we have the rules for natural numbers, which give us $\zero$, a successor operator $\succ(n)$ and an explicit recursor $\natrec(a,b,c)$ that incorporates the principle of complete mathematical induction

$$ \frac{}{\nat \val} \quad \frac{}{\zero \val} \quad \frac{}{\succ(\shouldsf{a}) \val} $$
$$\frac{\shouldsf{a \longmapsto a'}}{\natrec(\shouldsf{a,b,c}) \longmapsto \natrec(\shouldsf{a',b,c})}$$
$$ \frac{}{\natrec (\zero,\shouldsf{a,b}) \longmapsto \shouldsf{a}}$$ 
$$\frac{}{\natrec (\succ\shouldsf{(a),b,c}) \longmapsto \shouldsf{(c(a))}(\natrec(\shouldsf{a,b,c})).}$$

\vspace{4mm}

\noindent As it can be seen, the basic pattern here is that all type formers and their constructors are interpreted as values, eliminators preserve computation on their main argument, and eliminators transition to a particular information that was given to them for when their main argument is a certain constructor. 

The transition relation divides terms into two classes: those that are themselves values are called canonical, and those that are not themselves values but may eventually reach a value if we keep iterating the transition process are called non-canonical. This iteration is made precise with the notion of evaluation, a computation that halts when it finds the value that a term transitions to after an indefinite number of transition steps. I shall write $\shouldsf{a} \downarrow \shouldsf{a'}$ to mean that $\shouldsf{a}$ evaluates to $\shouldsf{a'}$. 
Since we are dealing with untyped computations, the evaluation of a term will not always terminate, but it suffices to characterize closed expressions as programs that when executed output the value they evaluate to. 
Notice that natural numbers are considered as canonical terms whereas, taking into account Frege's conception of numbers as extensions of concepts, one might expect such number terms to evaluate to lambda terms, the canonical terms representing value-ranges. This can be done with a technique known as Church encoding, with which we can represent the natural numbers using lambda notation. But the idea goes against the spirit of type theory in that it conflates two very distinct data types, numbers and functions. 

%-----------------------------------------------------------------------------------------------------------------------------------------
\subsection{Type and term equality} \label{typeeq}
%-----------------------------------------------------------------------------------------------------------------------------------------

Roughly, we explain what a type is by specifying the terms of that type following an idea that can be traced back to the constructive conception of set advocated by \cite{bishop1967foundations} and the interpretation of the intuitionistic logical constants proposed by \cite{heyting1934grundlagenforschung}. In fact, the first thing we do to explain what types are is distinguishing between canonical and non-canonical types. Any expression that evaluates to a canonical type is a type, and a canonical type is explained by prescribing what their canonical terms are. 

The meanings of the type and term equality judgments are mutually explained in a similar way. To account for the former, we assume that $A$ and $A'$ are closed expressions, and assert that

\begin{center}
	\textit{$A$ and $A'$ are equal types provided that they evaluate to the same canonical type up to type equality.}
\end{center}

\noindent This characterizes type equality as a relation between terms that behave in a certain expected way. The following rules illustrate how type equality can be explained for our canonical types: dependent function types, equality types, and natural number types

$$ \frac{A \downarrow \prod_{x : B} C \quad A' \downarrow \prod_{x : B'} C' \quad B \equiv B' \type \quad x : B \vdash C \equiv C'}{A \equiv A' \type} $$
$$ \frac{A \downarrow a =_B b \quad A' \downarrow a' =_{B'} b' \quad B \equiv B' \type \quad a \equiv a' : B \quad b \equiv b' : B }{A \equiv A' \type} $$
$$ \frac{A \downarrow \nat \quad A' \downarrow \nat}{A \equiv A' \type.} $$

\vspace{4mm}

\noindent While the above stipulations endow type equality with an intensional nature~\citep{dybjer2012program}, we often find an additional condition that determines a type uniquely by its terms

$$
\frac{a : A \vdash a : A' \quad a : A' \vdash a : A}{A \equiv A' \type.}
$$

Now, in order to give a full account of term equality we assume that we are given closed terms $a$ and $b$ as well as a type $A$, and evaluate them all, declaring that

\begin{center}
	\textit{$a$ and $b$ are equal terms of type $A$ provided that $a$ and $b$ evaluate to equal canonical terms of the canonical type which $A$ evaluates to.}
\end{center}

\noindent More concretely, term equality can be explained by laying down what it is to evaluate to equal canonical terms up to term equality for all the considered canonical types: 

$$ \frac{A \downarrow \prod_{x : B} C \quad a \downarrow \fun x M \quad a' \downarrow \fun x M' \quad x : B \vdash M \equiv M' : C}{a\equiv a' : A} $$
$$ \frac{A \downarrow b =_B b' \quad a \downarrow \refl(b) \quad a' \downarrow \refl(b') \quad b \equiv b' : B}{a \equiv a' : A} $$
$$ \frac{A \downarrow \nat \quad a \downarrow \zero \quad a' \downarrow \zero}{a\equiv a' : A}$$
$$\frac{A \downarrow \nat \quad a \downarrow \succ(M) \quad a' \downarrow \succ(M') \quad M\equiv M' : \nat}{a \equiv a': A.}$$
\vspace{4mm}

\noindent The semantics given above can be inductively extended to judgments occurring under a list of hypothesis, meaning that the validity of the hypothetical judgment 

$$ x : A  \vdash f : B $$

\noindent is subject to the validity of the non-hypothetical or categorical judgment $f[a/x] : B$ for every closed term $a : A$. In sum, the only way a hypothetical judgment can be given meaning to is by determining what categorical judgment they result in when all their free variables are replaced with closed terms. Frege's determination of the reference of functional expressions in \gga~\S\S29--31 is based on a similar idea, since the role that functional expressions play in the ideography is no different than that of open terms in type theory \citep{bentzen2020frege}. 

Because only closed terms are taken into account in the validation of a hypothetical judgment, it can be shown that the following principle, which can be regarded as a stronger version of the so-called Axiom K~\citep{streicher1993investigations}, obtains 

$$ x : A, p : x =_A x  \vdash p \equiv \refl(x) : x =_A x $$

\noindent because, for a closed term $a : A$, a closed equality term $b : a =_A a$ will always evaluate to a canonical term of $a =_A a$, which means that $p \downarrow \refl(x)$. Since $a =_A a$ is a canonical type, $c \equiv \refl(a) : a =_A a$ obtains. This principle is a straightforward consequence of the reflection rule
\begin{equation*}
\frac{\; p : x =_A y  }{x \equiv y : A \;}
%\frac{x : A, y : A  \vdash p : x =_A y  }{x : A, y : A \vdash x \equiv y : A}
\end{equation*}

\noindent which can be easily validated using a similar line of reasoning and expresses the fact that, semantically, every propositional equality is also a judgmental equality. 

%=========================================================================================================================================
\section{Sense and reference} \label{sense-and-reference}
%=========================================================================================================================================

It was \cite{dummett1978sense} who first noticed that Frege's theory of sense and reference could be rendered constructively by assuming that the sense of a singular term is related to its reference as a program is related to its value, and that the thesis that the reference of a sentence is a truth value should be entirely rejected, since, if we think of the sense of an expression not as a mode of presentation of its reference but as an effective method for determining it, that would imply a means of deciding whether a proposition is true or false.\footnote{
	\cite{moschovakis1994sense} adds to this interpretation a mathematical notion of recursive algorithm that allows for a more rigorous theory of sense. It can be seen as a forerunner of the theory of sense of reference of \cite{martinlof2001sense}, which is grounded in type theory instead.} 

Although the idea that a sentence refers to a truth value has no place in a constructive setting, \cite{martinlof2001sense} has observed that an interpretation of singular terms as terms and propositions as types in the sense of the meaning explanations allows for a more rigorous development of a theory of sense and reference for both singular terms and sentences. In this context, there is nothing better suited to  mediate the passage from the sense to the reference of an expression than the evaluation of a term to its canonical form. More concretely, the reference of a term $(\fun x \refl(x))(\zero)$ is the canonical term it evaluates to, $\refl(\zero)$. Programs are means for specifying their values and so are senses, since values are references. 
Expanding this idea to the sphere of types, Martin-L\"of notes that the reference of a sentence has to be a canonical proposition, a view that fits perfectly the meaning explanations, since an expression that is assigned to a type is interpreted as a sentence that expresses a proposition, and a type evaluates to a canonical type, which is taken as their referent. For example, both $(\fun x x)(\nat)$ and $\natrec(\zero,\nat,\fun x x(x))$ are equal types because they evaluate to the canonical type $\nat$. In other words, they are coreferential types but specify their value differently. 

If the passage from the sense to the reference of an expression is given through evaluation, how should one interpret the passage of an expression to its sense? I have suggested elsewhere that an expression $a$ comes to be known as a program when it is assigned to a type $A$.\footnote{
	See, for instance, \cite{bentzen2018constructive,bentzen202Xdifferent}.} 
The reason for this lies in an analogy with computer programs: a piece of code is only recognized as a program when it is correctly typechecked by a compiler of the programming language in question. If there is a single type mismatch error, such as an assignment of a value between two variables of different types, then we did not write a program after all. However, if we assume that a term $a$ has a sense only when $a : A$, we are at the same time excluding the possibility that terms lacking a reference have a sense. This has to do with the fact that $a : A$ obtains only when $a$ computes to a value according to our semantic stipulations, while not every term has a value. Consider a non-terminating term such as the infamous $\Omega$ term:

$$ (\fun x x(x))(\fun x x(x)).$$

\noindent Frege has always defended that in an exact science such as mathematics every expression must have a reference, but he never once entertained the idea that a sentence does not express a sense if it lacks a reference. Even in our setting, it is desirable to allow expressions to have a sense when they do not have a reference, for even a non-terminating term can be run just like any other program. 
I shall therefore say that an expression is recognizable as a program when it is interpretable by an operational semantics in the style of the previous section. That means that all that is needed for an expression to have a sense is a grammatical structure that conforms to the program syntax. Returning to our analogy with computer programs, we may say that a piece of code is a legitimate program when it passes a lexical check, when there are no missing required characters or unrecognizable tokens. Using Frege's semantic triangle, the tripartite distinction between expressions, programs, and values can be succinctly formulated as follows: 

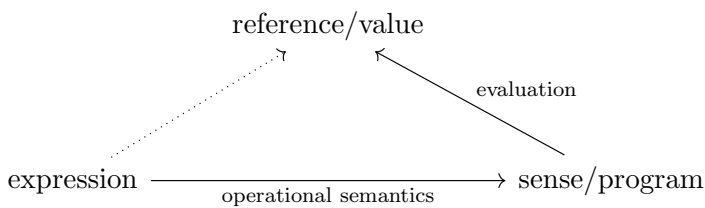
\begin{center}
	\begin{tikzcd}%[column sep=2.5em]
		& \text{reference/value} & \\
		&& \\
		\text{expression}  \arrow[uur,dotted,swap]\arrow[rr,swap,"\text{operational semantics}"] && \text{sense/program} \arrow[uul,swap,"\text{evaluation}"]
	\end{tikzcd}
\end{center}

Martin-L\"of is careful to distinguish between two different interpretations of equality of reference, claiming that a judgmental equality $a \equiv b : A$ says of the senses of $a$ and $b$ that they are coreferential and a propositional equality $a =_A b $, when proven to be true, says of the references of the expressions `$a$' and `$b$' that they are equal objects, because one is here suppressing the equality term $p$ that realizes the propositional equality via the assertion of $p : a =_A b $, the part of the judgment that he sees as not ``referentially transparent''. Having said that, I cannot see the grounds for this confusing distinction, considering that semantically every propositional equality is a judgmental equality. 
Since Martin-L\"of speaks of judgmental equality as an intensional relation, one may question if he is not simply leaving the semantics aside and speaking of the formalism itself, since, to give an example, while it is easy to construct a term of type $n + m =_\nat m + n$ via the principle of complete induction, there is no derivation of the judgmental equality $n + m \equiv m + n : \nat$ in the theory, since both sides of the equality do not evaluate to the same value, although they do for every particular instance of $n$ and $m$. Semantically, both propositional and judgmental equalities say that the senses of $a$ and $b$ have the same reference, because they amount to the fact that the values of programs $a$ and $b$ are one and the same. %There is no distinction,  

Martin-L\"of also claims that equality of sense is given by a relation of synonymy for expressions, which includes in particular renaming of bound variables and definitional stipulations such as $\pi = C / d$. I believe that this view of sameness of sense is essentially correct, but it lacks a theoretical motivation in the setting of the meaning explanations. The following computational explanation is able to provide a more compatible justification. Since to say that two expressions have the same sense is to say that two terms express the same program in type theory, we just have to justify an identity criterion for programs. Clearly, it seems appropriate to identify programs that have the same computational content. Hardly anyone would say that two programs are identical if one involves more computations than other, even though they may have the same execution value. With that in mind, I propose that $a$ and $b$ express the same program if they have the exact same computational behavior, if to evaluate $a$ is to evaluate $b$ and vice versa.  It is simply natural to identify two terms up to renaming of bound variables and definition unfolding, given that they are evaluated in the same way. 

%-----------------------------------------------------------------------------------------------------------------------------------------
\subsection{Equipollence}
%-----------------------------------------------------------------------------------------------------------------------------------------

It only remains to be asked how this account of sameness of sense can be linked to the equipollence principles that we have previously discussed. For the sake of argument, let us first assume that $A$ and $B$ are sentences, which is to say that they evaluate to canonical types. We have to show that $A$ and $B$ have the same sense when anyone who accepts $A$ will immediately accept $B$ and vice versa. One interpretation of the immediacy in the latter condition is that proposed in \cite{sundholm1994frege} 
\begin{equation*} 
\frac{a : A}{a : B} \quad\text{and}\quad \frac{b : B}{b : A}
\end{equation*}

\noindent that says that any construction that realizes $A$ is also a realizer for $B$ and vice versa. But that does not imply that $A$ and $B$ have the same computational behavior. Instead, it says that they have the same reference, since from this coextensionality it would follow that $A \equiv B \type$. That may seem odd from a Fregean perspective, but recall that sentences do not refer to truth values but to canonical types here.  

If we choose to interpret immediacy in a purely computational way instead, then we have that anyone who accepts $A$ will immediately accept $B$ if, without any reference to the notions of transition or evaluation, one is justified in an making assertion of the truth of $B$ from one's assertion of the truth of $A$, or, equivalently, if from one's construction of a term $a : A$ one can construct a (possibly new) term $b : B$ by means that do not involve the computation. It is plain that sameness of sense is given by equipollence under this view, for this condition describes precisely what it means for two types to have the same computational behavior. 
Finally, if we are given an open type 

$$x : A \vdash P(x) \type$$

\noindent then the corresponding principle for terms follows directly from the equipollence of types, because to say that two terms $a : A$ and $b : A$ are equipollent is to explain how anyone who accepts $P(a)$ will immediately accept $P(b)$, which we already did. Notice that only terms of a same type can be equipollent, like everything else in type theory. It does not make sense to ask whether, say, Julius Caesar is equipollent to $0$. 

If we were to interpret equality of reference for types in a strictly Fregean sense, that is, following the idea that two sentences are coreferential if they have the same truth value, then obviously the relation that we would be looking for would be that of logical equivalence, which, type theoretically, translates to the existence of two functions $f : A \to B$ and $g : B \to A$ between the types $A$ and $B$.\footnote{
	We typically write $A \to B$ as an abbreviation for $\Pi_{x : A} B$ when the type $B$ does not depend on $x : A$.} 
But it is worth stressing that actual equality of reference does not follow from logical equivalence, for two types may imply each other without computing to the same canonical type. Put differently, if we had a type universe $\mathcal{U}$, whose terms are smaller types, then the following principle of propositional extensionality 

$$ (A \to B) \to (B \to A) \to (A =_{\mathcal{U}} B) $$
would be false in the meaning explanations, just as it was in \bs~because $A$ and $B$ may have different contents. In general, we do not even have that $g(f(a))=_A a$ and $f(g(b))=_B b$ obtain, for every $a : A$  and $b : B$. The only thing logical equivalence tells us is that if we constructed a term of $A$ we know how to find another term of $B$ and vice versa. It comes with no other guarantees, as \cite{sundholm1994frege} notes.

%-----------------------------------------------------------------------------------------------------------------------------------------
\subsection{Equality and homotopy}
%-----------------------------------------------------------------------------------------------------------------------------------------

Over the last few years, homotopy type theory~\citepalias{hottbook} has emerged as a new foundation for mathematics that unites type theory and homotopy theory via a homotopical interpretation of type theory where a type $A$ is viewed as a space, a term $a : A$ as a point of the space $A$, an equality term $p : a =_A b$ as a path from point $a$ to point $b$ in the space $A$ and so on.\footnote{
	For a philosophical introduction to homotopy type theory, see \cite{ladymanpresnell2017hottfoundation} and \cite{bentzen2019what}.}

One of the main ingredients that reinforces this interpretation is the univalence axiom, which offers a formal justification for the common view among mathematicians that two mathematical objects are equal when they are isomorphic. Roughly, the univalence axiom implies that two types are identical just in case they are equivalent in a technical sense: 

$$ \mathsf{ua} : A \simeq B \to A =_{\mathcal{U}} B.$$
And this new sort of identification imposes a weaker understanding of propositional equality where in general a closed term $p : a =_A b$ does not entail $a \equiv b : A$, because univalence introduces a new canonical term $\mathsf{ua}(e)$ to the equality type that is not judgmentally equal to a reflexivity term, assuming that $e$ is an equivalence. Considering that $\refl(a)$ is interpreted as the constant path at the point $a$, univalence ensures the existence of non-trivial paths. 

Yet, this homotopical view of equality goes against the constructive theory of sense and reference that I have expounded, precisely because if propositional equality is to express coreferentiality, then the reference of a term cannot be its value in homotopy type theory. In the presence of non-trivial paths, that is, if we are allowed to construct a closed term $p : a = _A b$ such that $p \not\equiv \refl(a)$, then we have that $ a \downarrow a'$ and $b \downarrow b'$ but $a' \not\equiv b'$, meaning that the interpretation of propositional equality as equality of values is lost.

\bibliographystyle{plainnat}
{\linespread{0.7}\selectfont\bibliography{ref}}

\address

\end{document}